\documentclass[12pt]{article}%
\usepackage{graphicx}
\usepackage[intlimits]{amsmath}
\usepackage{latexsym}
\usepackage{amsfonts}
\usepackage{amssymb}%
 \makeatletter
\newfont{\footsc}{cmcsc10 at 8truept}
\newfont{\footbf}{cmbx10 at 8truept}
\newfont{\footrm}{cmr10 at 10truept}
\newtheorem{theorem}{Theorem}

\newtheorem{fact}[theorem]{Fact}

\newenvironment{proof}[1][Proof]{\noindent{\textbf {#1}  }}  {\hfill$\Box$\bigskip}

\begin{document}

\title{Spectral radius and Hamiltonicity of graphs}
\author{Miroslav Fiedler\thanks{Department of Computational Methods, Institute of
Computer Science, Academy of Sciences of the Czech Republic; \textit{e-mail:
fiedler@cs.cas.cz}} \ and Vladimir Nikiforov\thanks{Department of Mathematical
Sciences, University of Memphis, Memphis, TN, USA; \textit{e-mail:
vnikifrv@memphis.edu}}}
\maketitle

\begin{abstract}
Let $G$ be a graph of order $n$ and $\mu\left(  G\right)  $ be the largest
eigenvalue of its adjacency matrix. Let $\overline{G}$ be the complement of
$G.$

Write $K_{n-1}+v$ for the complete graph on $n-1$ vertices together with an
isolated vertex, and $K_{n-1}+e$ for the complete graph on $n-1$ vertices with
a pendent edge.

We show that:\medskip

If $\mu\left(  G\right)  \geq n-2,$ then $G$ contains a Hamiltonian path
unless $G=K_{n-1}+v;$ if strict inequality holds, then $G$ contains a
Hamiltonian cycle unless $G=K_{n-1}+e.$\medskip

If $\mu\left(  \overline{G}\right)  \leq\sqrt{n-1},$ then $G$ contains a
Hamiltonian path unless $G=K_{n-1}+v.$\medskip

If $\mu\left(  \overline{G}\right)  \leq\sqrt{n-2},$ then $G$ contains a
Hamiltonian cycle unless $G=K_{n-1}+e.$\medskip



\end{abstract}

\section{Introduction}

The spectral radius of a graph is the largest eigenvalue of its adjacency
matrix. In this note we give tight conditions on the spectral radius for the
existence of Hamiltonian paths and cycles. Other spectral conditions for
Hamiltonian cycles have been given in \cite{BuCh08}, \cite{Heu95},
\cite{KrSu03}, and \cite{Moh92}, but they all are far from our results.

In \cite{Ore60} Ore showed that if the inequality%
\begin{equation}
d\left(  u\right)  +d\left(  v\right)  \geq n-1 \label{or}%
\end{equation}
holds for every pair of nonadjacent vertices $u$ and $v,$ then $G$ contains a
Hamiltonian path. If the inequality (\ref{or}) is strict, then $G$ contains a
Hamiltonian cycle.

Write $K_{n-1}+v$ for the complete graph on $n-1$ vertices together with an
isolated vertex, and $K_{n-1}+e$ for the complete graph on $n-1$ vertices with
a pendent edge.

Note that $K_{n-1}+v$ has no Hamiltonian path as it is disconnected, and
$K_{n-1}+e$ has no Hamiltonian cycle as it has a vertex of degree $1.$ As it
turns out, these are the only graphs of order $n$ with such properties and
having maximum number of edges.

More precisely, Ore's approach can be used to prove the following extremal result.

\begin{fact}
\label{f1}Let $G$ be a graph with $n$ vertices and $m$ edges. If
\begin{equation}
m\geq\binom{n-1}{2} \label{siz}%
\end{equation}
then $G$ contains a Hamiltonian path unless $G=K_{n-1}+v.$ If the inequality
(\ref{siz}) is strict, then $G$ contains a Hamiltonian cycle unless
$G=K_{n-1}+e.$
\end{fact}

We can now easily deduce a straightforward spectral version of this assertion.

\begin{theorem}
\label{th1}Let $G$ be a graph of order $n$ and spectral radius $\mu\left(
G\right)  .$ If
\begin{equation}
\mu\left(  G\right)  \geq n-2, \label{mu}%
\end{equation}
then $G$ contains a Hamiltonian path unless $G=K_{n-1}+v.$ If the inequality
(\ref{mu}) is strict, then $G$ contains a Hamiltonian cycle unless
$G=K_{n-1}+e.$
\end{theorem}

\begin{proof}
Stanley's inequality \cite{Sta87}
\[
\mu\left(  G\right)  \leq-\frac{1}{2}+\sqrt{2m+\frac{1}{4}},
\]
together with (\ref{mu}), implies that%
\[
2m\geq\left(  n-\frac{3}{2}\right)  ^{2}-\frac{1}{4}=n^{2}-3n+2.
\]
Hence, we obtain
\[
m\geq\binom{n-1}{2}%
\]
with strict inequality if (\ref{mu}) is strict. Now Theorem \ref{th1} follows
from Fact \ref{f1}.
\end{proof}

The proof of Theorem \ref{th1} is so short because Stanley's inequality
becomes equality precisely for complete graphs together with isolated vertices.

Another, subtler condition for Hamiltonicity can be obtained using the
spectral radius of the complement of a graph.

\begin{theorem}
\label{th2}Let $G$ be a graph of order $n$ and $\mu\left(  \overline
{G}\right)  $ be the spectral radius of its complement. If
\begin{equation}
\mu\left(  \overline{G}\right)  \leq\sqrt{n-1}, \label{cmu1}%
\end{equation}
then $G$ contains a Hamiltonian path unless $G=K_{n-1}+v.$

If
\begin{equation}
\mu\left(  \overline{G}\right)  \leq\sqrt{n-2}, \label{cmu3}%
\end{equation}
then $G$ contains a Hamiltonian cycle unless $G=K_{n-1}+e.$
\end{theorem}

Our proof of Theorem \ref{th2} is based on the concept of $k$-closure of a
graph, used implicitly by Ore in \cite{Ore60}, and formally introduced by
Bondy and Chvatal in \cite{BoCh76}. We write $E\left(  G\right)  $ for the
edge set of a graph $G$ and $e\left(  G\right)  $ for $\left\vert E\left(
G\right)  \right\vert ;$ $d_{G}\left(  u\right)  $ stands for the degree of
the vertex $u$ in $G.$

Fix an integer $k\geq0.$ Given a graph, $G$ perform the following operation:
if there are two nonadjacent vertices $u$ and $v$ with $d_{G}\left(  u\right)
+d_{G}\left(  v\right)  \geq k,$ add the edge $uv$ to $E\left(  G\right)  .$ A
$k$\emph{-closure} of $G$ is a graph obtained from $G$ by successively
applying this operation as long as possible. Somewhat surprisingly, it turns
out that the $k$-closure of $G$ is unique, that is to say, it does not depend
on the order in which edges are added; see \cite{BoCh76} for details.

Write $\mathcal{G}_{k}\left(  G\right)  $ for the $k$-closure of $G$ and note
its main property:\medskip

$d_{\mathcal{G}_{k}\left(  G\right)  }\left(  u\right)  +d_{\mathcal{G}%
_{k}\left(  G\right)  }\left(  v\right)  \leq k-1$\emph{ for every pair of
nonadjacent vertices }$u$\emph{ and }$v$\emph{ of }$\mathcal{G}_{k}\left(
G\right)  .$\medskip

The usefulness of the closure concept is demonstrated by the following two
facts, due essentially to Ore \cite{Ore60}:

\begin{fact}
\label{f2} A graph $G$ has a Hamiltonian path if and only if $\mathcal{G}%
_{n-1}\left(  G\right)  $ has one.
\end{fact}

\begin{fact}
\label{f3} A graph $G$ has a Hamiltonian cycle if and only if $\mathcal{G}%
_{n}\left(  G\right)  $ has one.
\end{fact}

Armed with these facts we can carry out the proof of Theorem \ref{th2}%
.\bigskip

\begin{proof}
[Proof of Theorem \ref{th2}]For short, let $H=\mathcal{G}_{n-1}\left(
G\right)  .$ Assume that (\ref{cmu1}) holds but $G$ has no Hamiltonian path.
Then, by Fact \ref{f2}, $H$ has no Hamiltonian path either. Now the main
property of $\mathcal{G}_{n-1}\left(  G\right)  $ gives $d_{H}\left(
u\right)  +d_{H}\left(  v\right)  \leq n-2$\emph{ }for every pair of
nonadjacent vertices $u$ and $v$ of $H;$ thus,%
\[
d_{\overline{H}}\left(  u\right)  +d_{\overline{H}}\left(  v\right)
=n-1-d_{H}\left(  u\right)  +n-1-d_{H}\left(  v\right)  \geq n
\]
for every edge $uv\in E\left(  \overline{H}\right)  .$ Summing these
inequalities for all edges $uv\in E\left(  \overline{H}\right)  ,$ we obtain%
\[
\sum_{uv\in E\left(  \overline{H}\right)  }d_{\overline{H}}\left(  u\right)
+d_{\overline{H}}\left(  v\right)  \geq ne\left(  \overline{H}\right)  ,
\]
and since each term $d_{\overline{H}}\left(  u\right)  $ appears in the
left-hand sum precisely $d_{\overline{H}}\left(  u\right)  $ times, we see
that
\[
\sum_{v\in V\left(  \overline{H}\right)  }d_{\overline{H}}^{2}\left(
v\right)  =\sum_{uv\in E\left(  \overline{H}\right)  }d_{\overline{H}}\left(
u\right)  +d_{\overline{H}}\left(  v\right)  \geq ne\left(  \overline
{H}\right)  .
\]
Using the inequality of Hofmeister \cite{Hof88}, we obtain
\[
n\mu^{2}\left(  \overline{H}\right)  \geq\sum_{u\in V\left(  \overline
{H}\right)  }d_{\overline{H}}^{2}\left(  u\right)  \geq ne\left(  \overline
{H}\right)  .
\]
Since $\overline{H}\subset\overline{G},$ we have%
\[
\mu\left(  \overline{H}\right)  \leq\mu\left(  \overline{G}\right)  \leq
\sqrt{n-1},
\]
and so,
\[
n\left(  n-1\right)  \geq n\mu^{2}\left(  \overline{G}\right)  \geq n\mu
^{2}\left(  \overline{H}\right)  \geq ne\left(  \overline{H}\right)  .
\]
This easily gives $e\left(  \overline{H}\right)  \leq n-1$ and
\[
e\left(  H\right)  =\binom{n}{2}-e\left(  \overline{H}\right)  \geq\binom
{n-1}{2}.
\]

Since $H$ has no Hamiltonian path, Fact \ref{f1} implies that $H=K_{n-1}+v.$
If $G=H,$ the proof is completed, so assume that $G$ is a proper subgraph of
$K_{n-1}+v.$ Then $\overline{G}$ is a star $K_{1,n-1}$ of order $n$ together
with some additional edges; therefore $\overline{G}$ is connected. Hence, by
the Perron-Frobenius theorem,
\[
\mu\left(  \overline{G}\right)  >\mu\left(  K_{1,n-1}\right)  =\sqrt{n-1},
\]
contradicting (\ref{cmu1}) and completing the proof for Hamiltonian paths.

Assume now that (\ref{cmu3}) holds but $G$ has no Hamiltonian cycle. Using
Fact \ref{f3} and arguing as above, we see that
\[
e\left(  H\right)  >\binom{n-1}{2},
\]
and since $H$ has no Hamiltonian cycle, Fact \ref{f1} implies that
$H=K_{n-1}+e.$ If $G=H,$ the proof is completed, so assume that $G$ is a
proper subgraph of $K_{n-1}+e.$ Then $\overline{G}$ is a star $K_{1,n-2}$
together with some additional edges; therefore, $\overline{G}$ contains a
connected proper supergraph of $K_{1,n-2}.$ Hence, by the Perron-Frobenius
theorem,%
\[
\mu\left(  \overline{G}\right)  >\mu\left(  K_{1,n-2}\right)  =\sqrt{n-2},
\]
contradicting (\ref{cmu3}) and completing the proof.
\end{proof}

\textbf{Acknowledgement. }Part of this work was done during the meeting on
Spectra of Graphs in Rio de Janeiro, 2008. The authors are grateful to the
organizers and particularly to Prof. Nair de Abreu for the hospitality and the
wonderful atmosphere.

\end{document}